# A NOTE ON EDGE GUARDS IN ART GALLERIES

R. Nandakumar (nandacumar@gmail.com)

**Abstract:** We study the Art Gallery Problem with Edge Guards. We give an algorithm to arrange edge guards to guard only the inward side of the walls of any N-vertex simple polygonal gallery using at most roof (N/4) edge guards - a weakened version of Toussaint's conjecture on the number of edge guards that can guard an entire simple polygon.

Introductions to the full problem are in [1] and [2]). Here, we consider only simple polygonal (no holes, no self -intersections) galleries. Each polygon edge contains both its end vertices.

***Proposition:*** *For any N-vertex simple polygonal gallery displaying only paintings (only the inward side of the walls of the polygon need to be watched, not the entire interior region), we need at most roof (N/4) edge guards.*

**Edge visibility digraph:** A point P on the simple polygon is said to *see* another point Q on the polygon if no point of the line segment PQ is in the *exterior* of the polygon. For any simple polygon, we form a digraph as follows: the vertices of the digraph correspond to the edges of the polygon. If from points on a polygon edge E1, we can see the *whole* of another edge E2, the digraph has the directed edge from vertex corresponding to E1 to the vertex corresponding to E2. It is easy to see that this visibility between polygon edges *need not* be 2-way, hence the digraph.

*Note:* in the following, 'in-degree of an edge' of a polygon means the in-degree of the vertex corresponding to that edge in the edge visibility digraph of the polygon.

Every edge of any simple polygon is fully visible from both its neighboring edges, so the digraph always has a doubly linked circular list as its 'backbone' plus other edges.

**Observation 1:** In any simple polygon, no edge can have an in-degree of only 2. There are very few possible configurations of edges with in-degree 3 (call these **'weak edges'**). Figure 1 shows them. In-degree= 3 edges appear as pairs of adjacent edges (1A and 1B in figure 1) or sets of 4 adjacent edges (2 below, an **'arrowhead'**) but never as single edges or arbitrarily long strings of edges. In figure 1, we mark weak edges with slashes.

*Figure 1:*

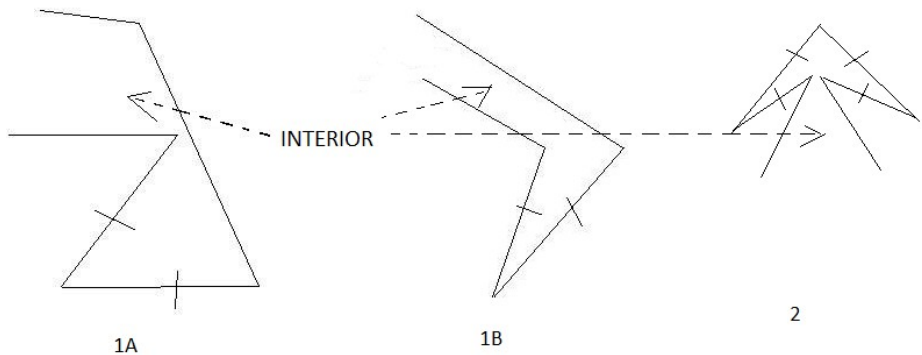

The above configurations (where 1A and 1B are not very different) do not generalize to arbitrarily long strings of weak edges. Indeed, such configurations as above apart, the other edges of the simple polygon have in-degree at least 4. Also note that configurations of weak edges are those regions where the edge visibility digraph is comparatively sparse and hence where more edge guards may be needed.

**Closest dominating vertex of an edge:** Consider an edge E (V1-V2) of a simple polygon. The half plane on the **inward** side of the infinite line containing this edge is potentially from where the inner wall of edge E can be seen. Of course, due to the presence of other edges, from most of this half plane, E is invisible or only partially visible. Among those vertices of the polygon (other than V1 and V2) lying in this half plane consider that vertex, say, V, whose shortest distance from edge E is the least. V is called the *closest dominating vertex* of E. Crucially, from V, the whole of E is necessarily visible.

**Observation 2:** Note all possible pairs of {edge, its closet dominating vertex} in the polygon. From this set, keep only those pairs where the closest dominating vertex has a concave angle. Consider all the line segments joining a chosen dominating vertex to the corresponding edge. We observe that no two of such line segments can have a proper intersection – although they may share an end point. In figure 2, V is the closest dominating vertex of E (the thick edge) and the shortest segment between them is shown dashed. No similar shortest segment can properly intersect it. Note: V could be the closest dominating vertex of an edge adjacent to E also (shown in figure 2 as a 'normal' thin edge) or even to other 'far' edges of the polygon.

*Figure 2*

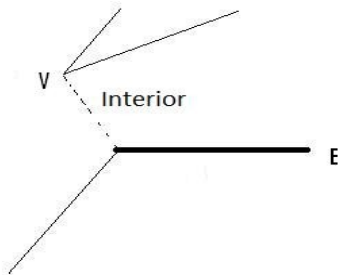

Each of these shortest line segments gives a **bottleneck** of the polygon that consists of an edge and the vertex that dominates it. Bottlenecks are of 2 types – (1) when the closest dominating vertex is a vertex of an adjacent edge of the dominated edge ('shallow bottlenecks') and (2) when the dominating vertex is *not* part of an adjacent edge of the edge it dominates ('normal bottlenecks').

**Observation 3:** Consider a simple polygon being partitioned into 2 sets of edges at any *normal* bottleneck. Call these edge sets the left and right sets. The minimum in-degree of any edge E1 (which may be 3 for groups of weak edges or 4 otherwise as noted above) is satisfied by edges in the same edge set as E1 plus those on the bottlenecks bounding this edge set. Indeed, (figure 3) if V and E form a bottleneck and if edge E1 on the left set is fully visible from an edge E2 on the right set, then, E is either fully visible from E and V or otherwise, there are edges in the left set which form a 'visibility channel' that passes thru the bottleneck and contains edge E. In the latter case, the edges forming the channel constitute another bottleneck and so E1 is not in an edge set bounded by E and V.

*Figure 3*

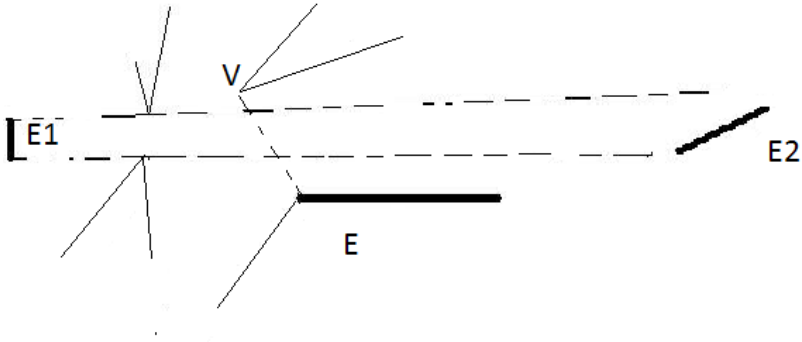

If any edge E of a polygon is made a guard, it automatically covers 3 edges, E itself and its two immediate neighbors. But, a guard at either of the two edges meeting at a dominating vertex in a normal bottleneck necessarily covers both its adjacent edges and also that edge dominated by the vertex, thus giving a minimum of 4 edges covered.

**Observation 4:** An edge guard at any given normal bottleneck usually covers the whole of *both* edge sets meeting there; if a vertex or edge is part of more than one bottleneck, a guard there usually covers several edge sets (far in excess of 4 edges). Exceptions are when there are shallow bottlenecks. We discuss them later (figure 6).

## Algorithm

Our approach to arranging edge guards to guard the boundary of the simple polygon is:
*Add edge guards incrementally such that after the r'th guard has been added, at least 4r edges of the polygon have been covered. If we can guarantee that this can be carried thru to the end, roof(N/4) guards will suffice for the full polygon (our main proposition).*

**Step 1:** Pre-process: identify and guard all the sets of edges with in-degree 3. We know the weak edge configurations need more guards so we deal with them first.

Figure 4 again shows the weak edge configurations. The weak edges are slashed; any of the darkened edges, if made a guard will cover all the weak edges in each case. Indeed, in each case, the guard covers at least 4 consecutive edges and often more than 4 edges. As each guard is added, we mark all edges it covers.

*Figure 4:*

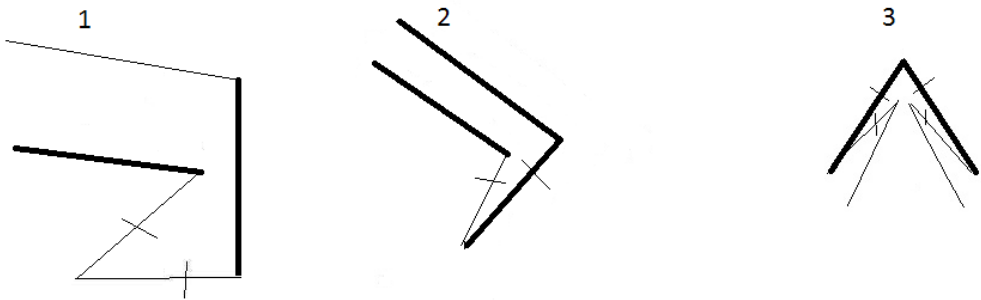

**Step 2:** Find normal bottlenecks and put guards (details follow).

At a normal bottleneck, let a guard be put at one of the edges meeting at the dominating (concave) vertex. This guard covers at least 4 edges of the polygon and separates edge sets with self-contained visibility – ie edges in each edge set do not depend on edges in the other edge set to achieve their minimum in-degree requirement.

At each nomal bottleneck there are 3 candidate edge guards possible – the 2 edges meeting at the dominating vertex and the dominated edge (experimentally, in most possible bottlenecks not 4 but at least 5 edges of the full polygon can be covered by one of these 3 possible edge guards). We also note that there is no need to guard all bottlenecks to cover all edges. Selecting from the available bottlenecks and from the possible guards at each bottleneck involves some heuristics described below. As each guard is put, we note how many edges get covered and mark them all. We try to ensure that each new guard covers 4 new edges. Note: The heuristics consider near-by bottlenecks but do not optimize globally.

*Figure 5:*

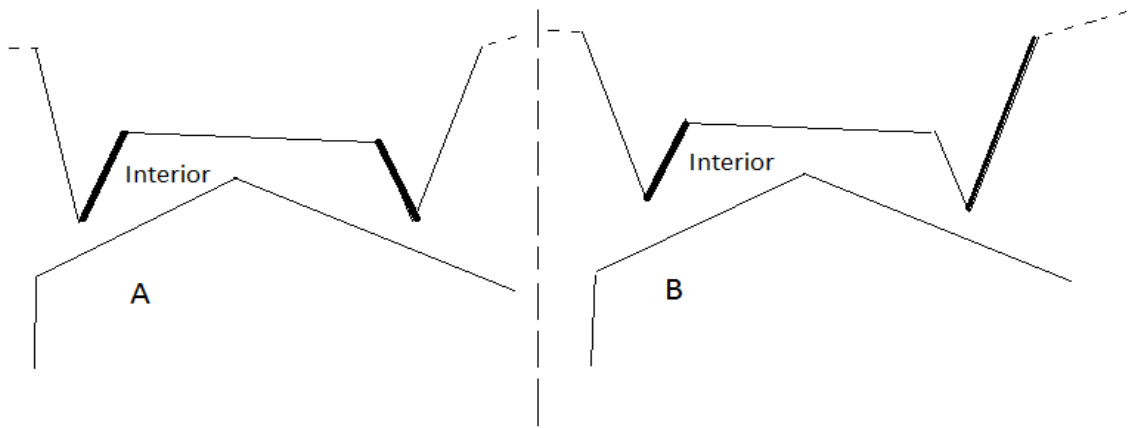

**Heuristics**: As shown in (say) fig 5, an edge may be part of two or more normal bottlenecks and a vertex may be part of several normal bottlenecks so a guard there covers several edge sets. Given a set of nearby normal bottlenecks, we try to put guards such that they cover distinct sets of edges (else, there is the risk: even if each guard does cover at least 4 edges, together they may not necessarily cover 8) and also such that together they do not leave out small sets of edges in between. Figure 5 shows a simple example. Clearly, from figure 5, arrangement B is to be chosen over A; and note that the bottleneck in the middle does not need a guard.

Note 1: Our approach is essentially local: we do consider nearby bottlenecks but not remote visibilities like those shown in figure 3. Indeed, we see that there is no need to optimize globally on the placement of guards - in every polygon we considered, there were multiple ways to cover all edges if roof(N/4) edge guards are available.

Note 2: As observation 2 says, the segments which mark bottlenecks of a polygon do not cut thru each other but may meet only at an end point. This ensures that the partition of the polygon into edge sets by normal bottlenecks is neat.

Now we consider *shallow bottlenecks with no weak edges* (these may occur within edge sets separated by normal bottlenecks). Figure 6A shows such a configuration. The vertex marked with a square dominates edge E1. It is easily seen that at least one of the edges meeting at that vertex (shown thick) covers a continuous string of 4 edges including E1. And since the in-degree of E1 is at least 4 (no edge is weak), there is at least one more edge in the part of the polygon which can cover E1 if guarded.

*Figure 6:*

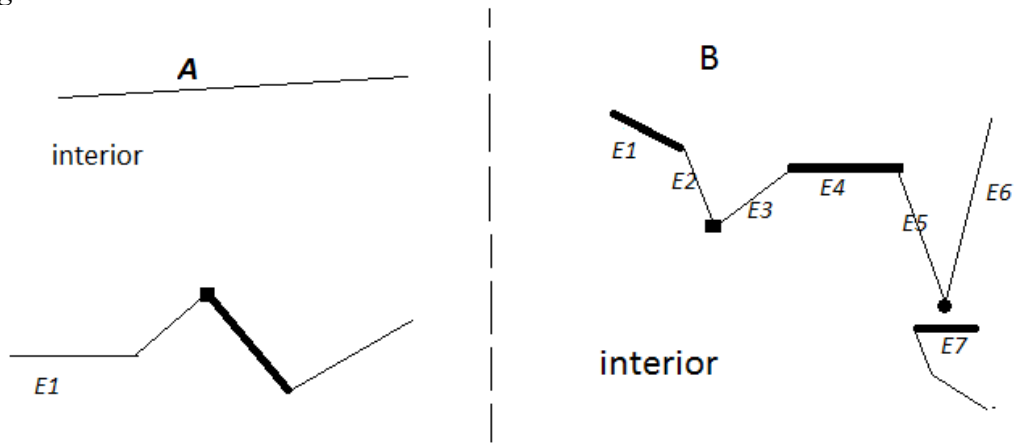

In 6B, the vertex marked with a dot formes a normal bottleneck with edge E7 and a large part of its left edge set is shown. In this edge set, the concave vertex marked with a small square dominates E4 (and possibly E1) in a shallow bottleneck. Due to this shallow bottleneck, none of the edges at the normal bottleneck at the dot can see the entire left edge set (E2 for instance is invisible). If there is more than one shallow bottleneck in an edge set (bounded by two normal bottlenecks), there could be parts of the edge set not covered even if we station guards at *both* normal bottlenecks bounding that edge set.

**Step 3:** If after step 2, any edges or edge sets remain uncovered (due to shallow bottlenecks), add guards for them.

After step 2 (inserting guards at normal bottlenecks), some portions of some edge sets might remained uncovered due to shallow bottlenecks. However, we observe the following: (1) such an edge set will have a large number of edges (otherwise with guards at its bounding bottlenecks, it will be fully covered) so addition of a guard or guards more can happen without exceeding the limit of a guard per 4 edges. (2) Any given edge set is a set of edges with high inter-connectivity in the visibility digraph (the in-degree of at least 4 is given to each edge by other edges in this edge set); so shallow bottlenecks lie in a well-interconnected part of the digraph and that helps. (3) Inserting edge guards at shallow bottlenecks on either side of an already guarded normal bottleneck could make the guard at the normal bottleneck unnecessary – a saving.

**Summing up:** At each weak edge group, an edge guard covering at least 4 successive edges can be put (Step 1). In step 2, each newly added normal bottleneck edge guard covers at least 4 edges - indeed, it covers at least 2 edge sets if there are no shallow bottlenecks; we also noted how shallow bottlenecks are not troublesome. Taken together, these observations imply that roof (N/4) edge guards suffice for the entire polygon.

**Remark:** If we bypass the preprocess stage of the algorithm (for in-degree = 3 edge groups) and only employ bottlenecks, there could be difficulties. Figure 6 shows a configuration of a bundle of arrows with narrow shafts joined at their bases. With the bottleneck-based approach, each unit arrow of 6 edges has two normal bottlenecks and we would need 2 edge guards per unit arrow (and hence more than roof(N/4) for the polygon) unless we make the bottleneck selection more sophisticated than discussed above. But with the preprocess on arrowheads, 1 guard per arrow will do and gives a total of much less than roof (N/4) guards.

*Figure 7:*

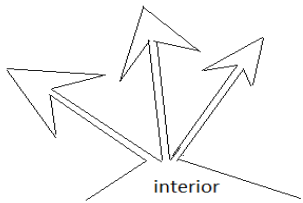

**Conclusion**: This three-stage algorithm (and heuristics) managed to cover all edges of any N vertex polygon we examined with roof (N/4) edge guards. A bonus: most polygons needed well below roof(N/4) guards. Indeed, for most polygons, the entire interior got covered along with the boundaries. As for those special polygons where only the boundary could be covered with less guards than the entire interior, our algorithm could cover the boundary with well below roof(N/4) guards and once the boundary was covered, a few additional guards were able to cover the entire interior - and the total number of guards was still within roof(N/4). So by attempting a weaker version of the edge guard conjecture, we might have got a handle on the full problem.